\pdfoutput=1
\documentclass[leqno,11pt]{amsart}
\usepackage[paperwidth=8.5in,paperheight=11in]{geometry} 
\usepackage{graphicx}
\usepackage[psdextra]{hyperref} 
\hypersetup{
   colorlinks=false,
   pdfborder={0 0 0},
}
\usepackage{tikz-cd}

\usepackage{enumitem}
\setlist{leftmargin=2pc}

\IfFileExists{mathrsfs.sty}{\usepackage{mathrsfs}}
{\IfFileExists{euscript.sty}{\usepackage[mathscr]{euscript}}
{\let\mathscr\mathcal}}

\theoremstyle{plain}
\newtheorem{theorem}{Theorem}
\newtheorem{corollary}[theorem]{Corollary}
\newtheorem{lemma}[theorem]{Lemma}
\newtheorem{proposition}[theorem]{Proposition}

\numberwithin{theorem}{section}

\theoremstyle{definition}

\newtheorem{example}[theorem]{Example}
 
\theoremstyle{remark}
\newtheorem{remark}[theorem]{Remark}

 \newcommand\Z{\mathbb {Z}}
\newcommand\Q{\mathbb {Q}} 
\newcommand\C{\mathbb {C}} \newcommand\A{\mathbb {A}}
\renewcommand\P{\mathbb {P}} 

\newcommand{\rk}{\operatorname{rank\,}_{\Z}}

\newcommand{\Gal}{\operatorname{Gal}}

\newcommand{\rank}{\operatorname{rank}}

\newcommand{\Qbar}{\overline{\Q}}

\newcommand{\ph}{\varphi}
\def\<{\langle}
\renewcommand{\>}{\rangle}

\newcommand{\Tr}{\operatorname{Tr}}

\newcommand{\Mod}[1]{\ (\mathrm{mod}\ #1)}

\numberwithin{equation}{section}

\title[Sextic Dirichlet Twists]%
{Ranks of Elliptic Curves in Cyclic Sextic Extensions of~$\Q$}
\date{}

\author[H.~Kisilevsky]{Hershy Kisilevsky}
\address[H.~Kisilevsky]%
{Department of Mathematics and Statistics and CICMA\\
     Concordia University \\
     1455 de Maisonneuve  Blvd. West\\
     Montr\'eal, Quebec, H3G 1M8, CANADA}
\email{kisilev@mathstat.concordia.ca}
\author[M.~Kuwata]{Masato Kuwata}
\address[M.~Kuwata]%
{Faculty of Economics \\
Chuo University\\
Hachioji-shi, Tokyo 192-0393,  Japan}
\email{kuwata@tamacc.chuo-u.ac.jp}

\subjclass[2020]{11G05, 14G05, 14G25, 11G40, 14J28}
\date{October 31, 2023}

\begin{document}

\maketitle

\begin{abstract}
For an elliptic curve $E/\Q$ we show that there are infinitely many cyclic sextic extensions $K/\Q$ such that the Mordell-Weil group $E(K)$ has
rank greater than the subgroup of $E(K)$ generated by all the $E(F)$ for the proper subfields $F \subset K$. For certain curves $E/\Q$ we show that the number of such fields $K$ of conductor less than $X$ is $\gg\sqrt X$.
\end{abstract}

\section{Introduction}    
Let  $E/\Q$ be an elliptic curve defined over $\Q$ with $L$-function $L(E,s).$ In \cite{FKK}, we examined the rank of the Mordell-Weil groups $E(K)$ as $K$ ranged over cyclic cubic extensions of $\Q.$ In view of the Birch \& Swinnerton-Dyer conjecture this is equivalent to considering the order of vanishing of the $L$-function
$L(E/K,s)$ at $s=1.$ More generally, if $\chi$ is a primitive Dirichlet character of order $d$,  $K= K_{\chi}$ is the cyclic extension of $\Q$ of degree $d$ corresponding to $\chi$, and $L(E/\Q,s,\chi)$
is the twist of $L(E,s)$ by $\chi$, then (see \cite{Ro}), the order of vanishing of $L(E/\Q,s,\chi)$ at $s=1$ is conjectured to be the rank of the
``$\chi$-component'' $E(K)^{\chi}$ of $E(K)$.  Here $\rk E(K)^{\chi}= \dim_{\C} \bigl(\C\otimes E(K)\bigr)^{\chi}$ is the dimension of the $\chi$ eigenspace of $\C\otimes_{\Z} E(K)$ as a $\Gal (\Qbar/\Q)$-space.
Kato's important result \cite{Kato} generalizing Kolyvagin's theorem \cite{Ko} shows
that if the $\chi$-component of $E(K_{\chi})$ has positive rank, then
$L(E/\Q,1,\chi)=0 $ (see Scholl \cite{Sch}) . 

In this article, we consider the case that $d=6$, and show that for any elliptic curve $E/\Q$, there are infinitely many primitive sextic characters
(and hence infinitely many cyclic sextic extensions $K=K_{\chi}$ of $\Q$) for which  $\chi$-component of $E(K_{\chi})$ has positive rank and so $L(E/\Q,1,\chi)=0. $
For certain curves $E/\Q$ we show that the number of such fields $K$ of conductor less than $X$ is $\gg\sqrt X$ which up to logarithmic factors is  the conjectured frequency of vanishing predicted in \cite{D-F-K}, \cite{D-F-K-07}, \cite{MR}, and \cite{Berg-Ryan-Young}.
 
 The method of proof is algebraic and is adapted from the arguments in \cite{FKK}.

\subsection*{Acknowledgements}
Kuwata was supported by JSPS KAKENHI Grant Numbers JP19K03427, and by the Chuo University Overseas Research Program.  
Part of this work was done while Kuwata was visiting Boston University.  
We thank the anonymous referee for many useful comments and suggestions.

\section{Vanishing Sextic Twists and Rational Points of a $K3$ surface}

In order to find points on $E/\Q$ defined over some cyclic 
extension $K_{\chi}/\Q$, we will
define an auxiliary variety of higher dimension whose $\Q$-rational points
correspond to points on $E$ defined over some cyclic extensions of $\Q$. 

For the case $d=6$, consider the automorphism $\rho$ of $E\times E$ defined by
\[
\rho: E\times E \to E\times E; \  (P,Q)\mapsto (Q,Q-P).
\]
Clearly, $\rho$ is defined over $\Q$, and it is readily verified that the order of $\rho$ is~$6$.

\begin{lemma}\label{lem:stabilizer}
The points of $E\times E$ over $\Qbar$ with nontrivial stabilizer are as follows.
\begin{enumerate}[label=\textup{(\roman*)}]
\item $\{(O,O)\}$; Stabilizer $=\<\rho\>$.
\item $\{(T_3,-T_3)\mid T_3\in E[3](\Qbar)\}\setminus\{(O,O)\}$; Stabilizer $=\<\rho^{2}\>$.
\item $\{(T_2,T'_2)\mid T_2,T'_2\in E[2](\Qbar)\}\setminus\{(O,O)\}$; Stabilizer $=\<\rho^{3}\>$.
\end{enumerate}
\end{lemma}

\begin{proof}
Let $P$ and $Q$ be points on $E(\Qbar)$, and suppose that $\rho^{i}(P,Q)=(P,Q)$ for some $i$, $1\le i\le 5$.  In other words, $(P,Q)$ equals one of the following five pairs:
\[
\quad (Q,Q-P), \quad (Q-P,-P), \quad 
(-P,-Q), \quad (-Q,P-Q), \quad\text{or} \quad (P-Q,P).
\]
If $(P,Q)=(Q,Q-P)$ or $(P-Q,P)$, then we have $P=Q=O$, which is case (i).  If $(P,Q)=(Q-P,-P)$ or $(-Q,P-Q)$, then we have $[3]P=O$ and $Q=-P$, which is case (ii).  If $(P,Q)=(-P,-Q)$, then we have $[2]P=[2]Q=O$, which is case (iii).
\end{proof}

Let $S_{6}(E)=E\times E/\<\rho\>$ be the quotient surface defined over~$\Q$, we denote by $[P,Q]$ the image in $S_{6}(E)$ of the point $(P,Q) \in E\times E$.  The group $\<\rho\>$ acts freely away from the set $\Sigma$ consisting of twenty four points in Lemma~\ref{lem:stabilizer}.  We denote by $S_{6}(E)^{\circ}$ the open subset of $S_{6}(E)$ obtained by removing the image of $\Sigma$.  Since $\Sigma$ is collectively defined over~$\Q$, $S_{6}(E)^{\circ}$ is also defined over~$\Q$.  A point $[P,Q]\in S_{6}(E)^{\circ}$ is defined over $\Q$ if and only if it is fixed by any element $\sigma$ in $\Gal(\Qbar/\Q)$.  In other words, the action of $\sigma$ is compatible with that of $\rho$, i.e., $\sigma(P,Q)$ equals $\rho^{i}(P,Q)$ for some $0\leq i \le5$.  More precisely, we have the following.

\begin{lemma}\label{lem:rat-pts-s6}
If a point $[P,Q]$ belongs to the set of rational points
$S_{6}(E)^{\circ}(\Q)$, then $(P,Q)$ satisfies one and only one of the following{\/\upshape:}
\begin{enumerate}[label=\textup{(\roman*)}]
\item 
$P$ and $Q$ are defined over $\Q$.
\item 
$P$ and $Q$ are both defined over some quadratic extension $\Q(\sqrt{\delta})/\Q$.  If $\tau\in\Gal(\Q(\sqrt{\delta})/\Q)$ is the generator, then 
we have $\tau(P)=-P$ and $\tau(Q)=-Q$.
\item 
$P$ and $Q$ are both defined over some cyclic cubic extension 
$K_{3}/\Q$. If we choose a suitable generator $\sigma_{3}\in\Gal(K_{3}/\Q)$, then 
we have $\sigma_{3}(P)=-Q$ and $\sigma_{3}(Q)=Q-P$.
\item 
$P$ and $Q$ are both defined over some cyclic sextic extension 
$K_{6}/\Q$. If we choose a suitable generator $\sigma_{6}\in\Gal(K_{6}/\Q)$, then 
we have $\sigma_{6}(P)=Q$ and $\sigma_{6}(Q)=Q-P$.
\end{enumerate}
\end{lemma}

\begin{proof}
Let $P$ and $Q$ be points on $E(\Qbar)$, and suppose that $[P,Q]\in
S_{6}^{\circ}(\Qbar)$.  If $\sigma\in \Gal(\Qbar/\Q)$, then
$\sigma(P,Q)=(\sigma(P),\sigma(Q))$ equals $\rho^{i}(P,Q)$ for some $i$, $0\le i\le 5$.  Note that these six pairs are distinct since $[P,Q]$ is in
$S_{6}^{\circ}(\Q)$.  We can therefore define a map $\psi:\Gal(\Qbar/\Q)\to
\<\rho\>$.  Since the automorphism $\rho$ is defined over $\Q$, it
commutes with any element of $\Gal(\Qbar/\Q)$.  Thus, if
$\sigma_{1}(P,Q)=\rho^{i}(P,Q)$ and $\sigma_{2}(P,Q)=\rho^{j}(P,Q)$, then
$(\sigma_{1}\circ\sigma_{2})(P,Q)=\sigma_{1}\bigl(\sigma_{2}(P,Q)\bigr) =\sigma_{1}\bigl(\rho^{j}(P,Q)\bigr)
=\rho^{i}\bigl(\rho^{j}(P,Q)\bigr) =\rho^{i+j}(P,Q)$.  This shows that the map 
$\psi$ is a homomorphism.

Let $K$ be the Galois extension of $k$ corresponding to $\ker\psi$ via 
Galois theory.  Then $\Gal(K/\Q)$ is isomorphic to a subgroup of 
$\<\rho\>$, which is a cyclic group of order~$6$.  

If $\Gal(K/\Q)=\{\text{id}\}$, then $K=\Q$, and both $P$ 
and $Q$ are defined over $k$.  This is the case~(i).

If $\Gal(K/\Q)\simeq\<\rho^{3}\>$, then $K$ is a quadratic extension of 
$\Q$ since $\rho^{3}$ is of order~$2$.  Let $\tau\in\Gal(\Qbar/\Q)$ be an element whose image in $\Gal(K/\Q)$ generates $\Gal(K/\Q)$.  Then $\psi(\tau)=\<\rho^{3}\>$.  This shows  
that $\tau(P,Q)=\rho^{3}(P,Q)=(-P,-Q)$.  This is the case~(ii).

If $\Gal(K/\Q)\simeq\<\rho^{2}\>$, then $K$ is a cyclic cubic extension of 
$\Q$ since $\rho^{2}$ is of order~$3$.  Let $\sigma_{3}\in\Gal(\Qbar/\Q)$ be an element whose image in $\Gal(K/\Q)$ generates $\Gal(K/\Q)$.  Then $\psi(\sigma_{3})=\<\rho^{2}\>$.  This shows  
that $\sigma_{3}(P,Q)=\rho^{2}(P,Q)=(Q-P,-P)$.  This is the case~(iii).

If $\Gal(K/k)\simeq\<\rho\>$, then $K$ is a cyclic sextic extension.  
Let $\sigma_{6}\in\Gal(\Qbar/\Q)$ be an element whose image in $\Gal(K/\Q)$ 
maps to $\rho$ by $\psi$.  Then we have $\sigma_{6}(P,Q)=(Q,Q-P)$.  This implies that $\sigma_{6}(P)=Q$ and $\sigma_{6}(Q)=Q-P$ which is the case~(iv).
\end{proof}

\begin{lemma}\label{lem:s6-K3}
The minimal nonsingular model $\widetilde{S}_{6}(E)$ of $S_{6}(E)=E\times E/\<\rho\>$ is a K3 surface, whose Picard number is either $20$ or $19$ depending on whether $E$ has complex multiplication or not.
\end{lemma}

\begin{proof}
Each of the points of nontrivial stabilizer in Lemma~\ref{lem:stabilizer} becomes a rational double point (or a simple surface singularity) in the quotient $S_{6}(E)$, and the number and type of singularities are as follows.
\begin{enumerate}[label=(\roman*)]
\item Type $A_{5}$ at one point $[O,O]$.
\item Type $A_{2}$ at four points of the form $[T_{3},2T_{3}]$, $T_{3}\in E[3]\setminus\{O\}$.
\item Type $A_{1}$ at five points of the form $[T_2,T'_2]$, $T_2,T'_2\in E[2]$, $(T_2,T'_2)\neq (O,O)$.
\end{enumerate}
The minimal resolution $\widetilde{S}_{6}(E)$ is obtained by blowing-up these singular points and replacing them with a tree of smooth rational curves, with intersection pattern dual to a Dynkin diagram of A-D-E singularity type.  Since the Galois conjugate points have the same type of singularities, the minimal resolution can be done over~$\Q$, and $\widetilde{S}_{6}(E)$ is defined over~$\Q$.  

It is readily verified that $S_{6}(E)$ satisfies ``Condition (K)'' in \cite[Th.~2.4]{katsura}; clearly $\<\rho\>$ is symplectic, has no fixed curve, and has only isolated fixed points which are rational double points.  So, the surface $\widetilde{S}_{6}(E)$ is a $K3$ surface.

As a result of the resolution of singularities, we obtain $1\times 5+4\times 2+5\times 1=18$ exceptional divisors.  If $E$ does not have complex multiplication, the N\'eron-Severi group of $E\times E$ is generated by the classes of divisors $\{O\}\times E$, $E\times \{O\}$, and $E\times E$.  Since the action of $\rho$ moves these classes transitively, the image of these divisors in the N\'eron-Severi group of $\widetilde{S}_{6}(E)$ generate a free abelian group of rank~$1$.  Thus, in this case the Picard number of $\widetilde{S}_{6}(E)$ equals $18+1=19$.  If $E$ has complex multiplication, then there is a nontrivial endomorphism $\ph:E\to E$, and the orbit of the divisor $E\times \ph(E)$ gives an extra independent element of the N\'eron-Severi group.  Thus, the Picard number of $\widetilde{S}_{6}(E)$ equals~$20$ in this case.
\end{proof}

Next, we fix a model of $E$ in $\P^{2}$ and see how $\rho$ acts on a pair of points $(P,Q)$ geometrically.  To do so, choose a model of $E$ over $\Q$ that is a cubic curve in $\P^{2}$ such that the origin $O\in E(\Q)$ is an inflection point.  For example, we may take a Weierstrass model of $E$.  Then, we have the property that three points $P_{1}$, $P_{2}$, and $P_{3}$ are collinear if and only if $P_{1}+P_{2}+P_{3}=O$.  Also, the inversion map $[-1]:P\mapsto -P$ on $E$ can be extended to an involution $\iota:\P^{2}\to \P^{2}$.  

Take two points $P$, $Q\in E(\Qbar)$, and let $R=Q-P$.  Then $P$, $-Q$, and $R$ are collinear, so we denote by $L$ the line passing through these three points.  This also implies that the line $\iota(L)$ passes through the three points $-P$, $Q$, and $-R$.  This implies that the pair $(P,Q)$ determines a degenerate conic  $C_{L}=L\sqcup \iota(L)$ consisting of a pair of conjugate lines.   (See Figure~1.)  Clearly, any pair of points $\rho^{i}(P,Q)$, $0\le i\le 5$ determines the same degenerate conic $C_{L}$.  Note that any pair $(Q,P)$ also determines the same $C_{L}$.

\begin{figure}[h]\caption{\null}
\centering
\includegraphics[scale=0.6]{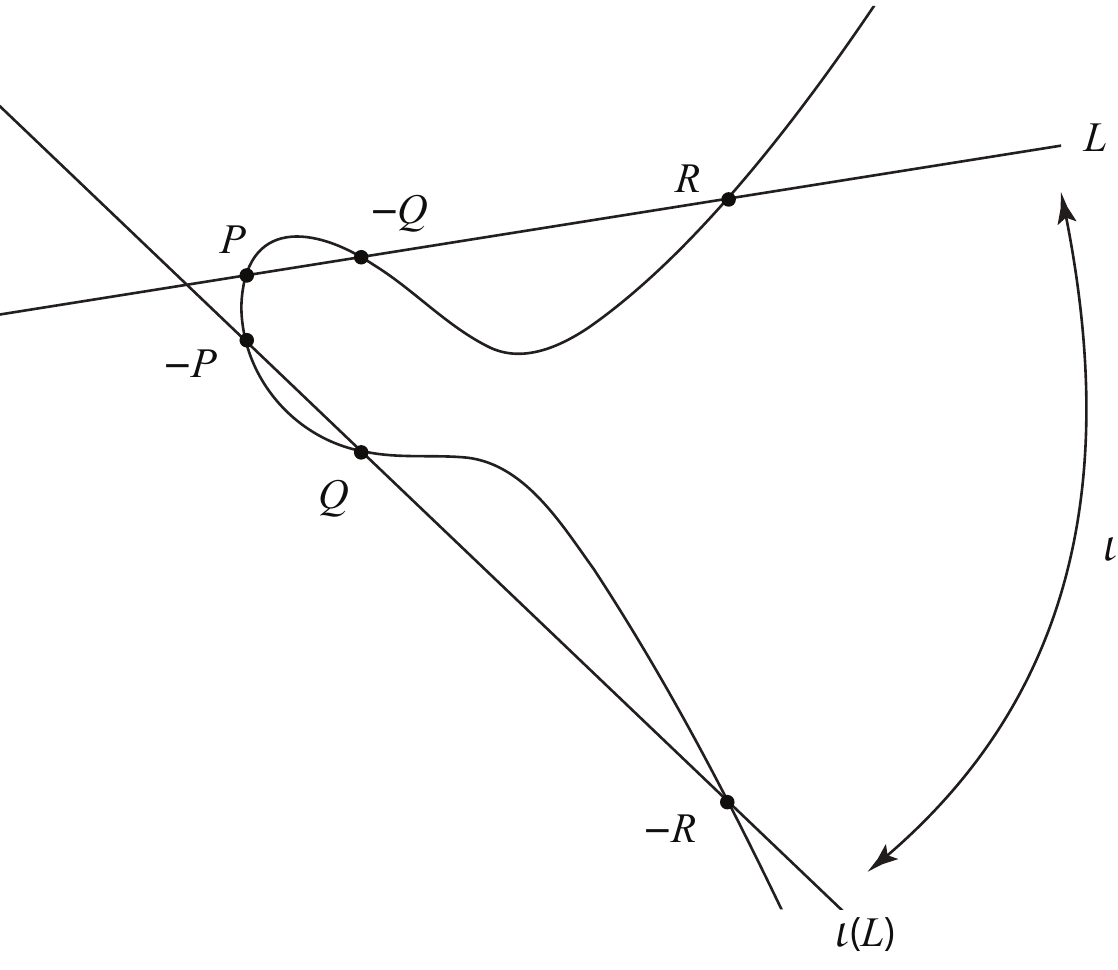}
\end{figure}

Conversely, take a line $L$ in $\P^{2}$ and let $C_{L}=L\sqcup \iota(L)$ be the degenerate conic.  Suppose $L$ is not the tangent line at a $2$-torsion point or a $3$-torsion point of $E$, and choose a point $P\in E\cap L$ and a point $Q\in E\cap \iota(L)$ different from $-P$.  Then, $(P,Q)$ determines $C_{L}$ in the above sense.
%
%

Now, take a degenerate conic $C_{L}=L\sqcup \iota(L)$ that is defined over~$\Q$.  Then, each line $L$ and $\iota(L)$ is a priori defined over a quadratic extension.  Suppose $L$ is indeed defined over a quadratic extension $K_{L}/\Q$ and not over $\Q$ itself, and let $\tau\in\Gal(\Qbar/\Q)$ be an element that generates $\Gal(K_{L}/\Q)=\Gal(\Qbar/\Q) \big/\!\Gal(\Qbar/K_{L})$.  Then, $\tau$ interchanges $L$ and $\iota(L)$, in other words, $\tau$ coincides with~$\iota$ on the points of $C_{L}$.
Choose $P\in E\cap L$ and $Q\in E\cap \iota(L)$ such that $Q\neq-P$, and let $R=Q-P$.  The intersection $E\cap C_{L}$ may be regarded as an effective divisor of degree~$6$ defined over~$\Q$, and may be written 
\[
E\cap C_{L}=(P)+(-P)+(Q)+(-Q)+(R)+(-R) \quad \text{as divisors}.
\]
Here, ``$(P)$'' means a divisor, and ``$+$'' means the formal sum, not in the sense of the group operation of $E$. 
Since each of the divisors $D_{P}=(P)+(-P)$, $D_{Q}=(Q)+(-Q)$, and $D_{R}=(R)+(-R)$ is invariant under~$\tau$, and the sum $D_{P}+D_{Q}+D_{R}$ is defined over~$\Q$, they are defined over some extension of $\Q$ whose Galois group is a subgroup of symmetric group $\mathfrak{S}_{3}$.

\begin{lemma}\label{lem:cyclic-sextic}
Suppose the field of definition of $L$ is a quadratic extension $K_{L}/\Q$, and the field of definition of the divisors $D_{P}$, $D_{Q}$ and $D_{R}$ is  a cyclic cubic extension $K_{3}/\Q$.  Then the field of definition $K_{6}$ of $P$ is a quadratic extension of $K_{3}$ such that $K_{6}/\Q$ is a cyclic sextic extension.  We may choose a generator $\varrho\in \Gal(K_{6}/\Q)$ such that 
\[
\begin{array}{lll}
\varrho(P)=Q, & \varrho^{2}(P)=R, & \varrho^{3}(P)=-P,  \\[\medskipamount]
\varrho^{4}(P)=-Q,\quad & \varrho^{5}(P)=-R, \quad& \varrho^{6}(P)=P.
\end{array}
\]
\end{lemma}

\begin{remark}
The action of $\varrho$ is indicated in Figure~2.
\end{remark}

\begin{figure}[h]\caption{\null}
\centering
\includegraphics[scale=0.64]{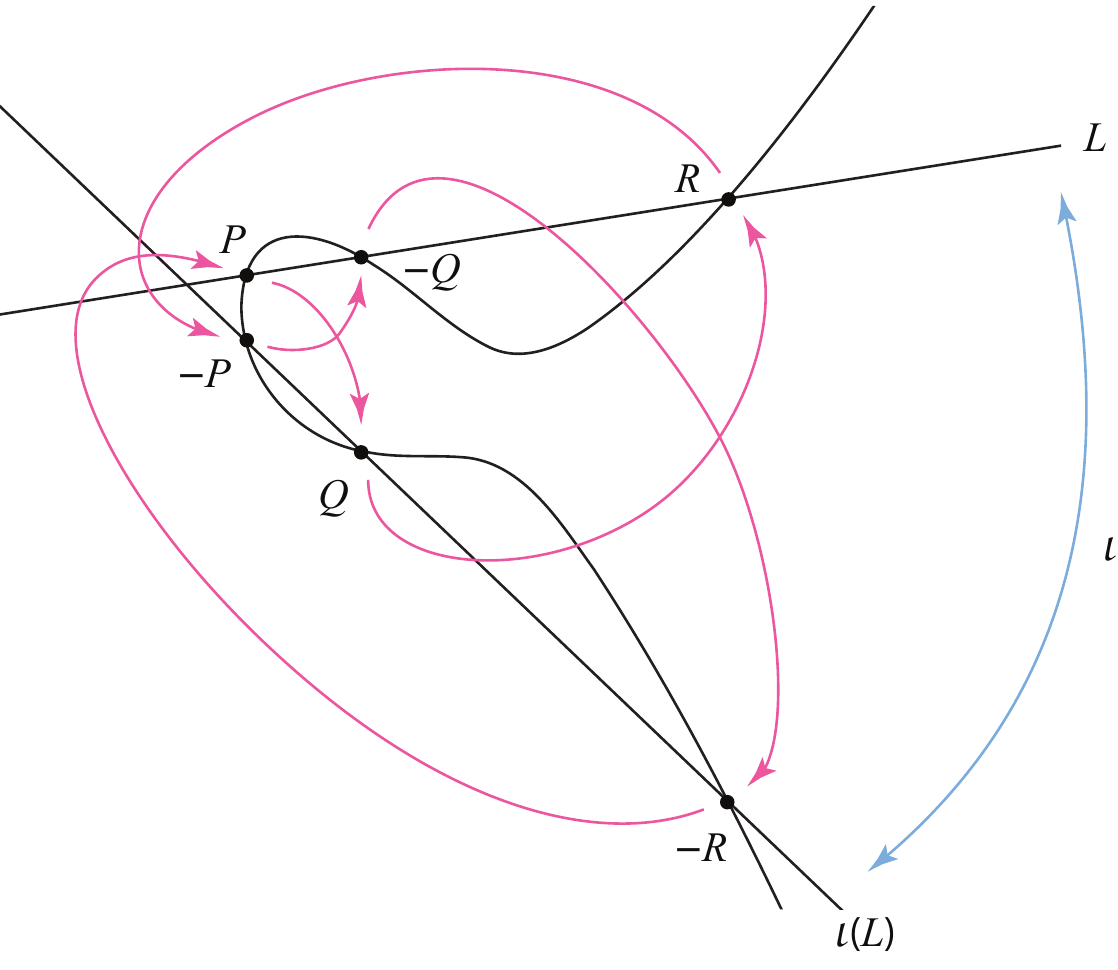}
\end{figure}

\begin{proof}
Let $\sigma\in\Gal(\Qbar/\Q)$ be an element that generates $\Gal(K_{3}/\Q)$. 
Since the Galois action commutes with the group operation of~$E$, $\sigma$ commutes with $\iota$.  But $\tau\in\Gal(\Qbar/\Q)$ coincides with $\iota$ on the points of $C_{L}$, and thus $\sigma$ commutes with $\tau$.  Since the Galois group of the Galois closure of $K_{L}K_{3}$ is generated by $\sigma$ and $\tau$, the Galois group $\<\sigma,\tau\>$ is a cyclic group of order~$6$.  So, $K_{6}=K_{L}K_{3}$ is a cyclic sextic extension of~$\Q$.

Let $\varrho\in\Gal(K_{6}/\Q)$ be a generator.  Then, since $\varrho^{3}=\iota$, $\varrho$ moves $L$ to $\iota(L)$.  If $\varrho(P)=-P$, then $D_{P}$ is fixed by $\varrho$, which is a contradiction.  Thus, $\varrho(P)$ equals either $Q$ or $-R$.  If $\varrho(P)=-R$, we replace $\varrho$ by $\varrho^{5}$, and we may assume $\varrho(P)=Q$.  Then, by the same reason, $\varrho(Q)$ must be on the line $L$, and it cannot be equal to neither $P$ nor $\iota(Q)$.  This implies $\varrho(Q)=R$, and $\varrho^{k}(P)$, $k=3,4,5,6$ are as in the statement since $\varrho^{3}=\iota$.
\end{proof}

\begin{remark}
The assumption of Lemma~\ref{lem:cyclic-sextic} implies that points $\pm P$, $\pm Q$, and $\pm R$ are all distinct.
\end{remark}

\begin{lemma}\label{lem:rankE(K6)}
Suppose $[P,Q]$ is a rational point in $S_{6}(E)^{\circ}(\Q)$ such that the field of definition of $P$ is a cyclic sextic extension $K_{6}/\Q$. If $P \in E(K_6)$ has infinite order, then $P$ is not in the subgroup $E(\Q)+ E(K_2)+ E(K_3)$, where $K_{2}$ and $K_{3}$ are unique quadratic and cubic subfields of~$K_{6}$.
\end{lemma}

\begin{proof}
By assumption, we are in the case (iv) of Lemma~\ref{lem:rat-pts-s6}, and we may choose a generator $\varrho\in \Gal(K_{6}/\Q)$ such that $\varrho(P)=Q$ and $\varrho(Q)=Q-P$.  In particular, we have $\Tr_{\varrho}(P)=O$, where $\Tr_{\varrho}(P)=\sum_{k=0}^{5}\varrho^{k}(P)$ denotes the trace of $P$ by $\varrho$.  

Suppose $P$ is in $E(\Q)+ E(K_2)+ E(K_3)$, and write $P=P_{1}+P_{2}+P_{3}$ with $P_{i} \in E(K_{i})$ for $i=1,2,3$, where $K_{1}=\Q$.  Since $K_{i}$ is the fixed field of $\varrho^{i}$ for $i=1,2,3$, we have $\varrho^{2}(P_{2})=P_{2}$ and $\varrho^{3}(P_{3})=P_{3}$.  Define
\[
P'_{1}=[6]P_{1} + \Tr_{\varrho}(P_{2}) + \Tr_{\varrho}(P_{3}), \quad
P'_{2}=[6]P_{2}-\Tr_{\varrho}(P_{2}), \ \text{ and } \
P'_{3}=[6]P_{3}-\Tr_{\varrho}(P_{3}),
\]
and let $P'=[6]P$.  Then, we have $P'=P'_{1}+P'_{2}+P'_{3}$ with $P'_i \in E(K_i)$, $i=1,2,3$.  Moreover, we have $\varrho(P'_1)=P'_1$, $\varrho(P'_2)=-P'_2$, and $P'_3+\varrho(P'_3)+\varrho^2(P'_3)=O$.  

We have, on the one hand, $\Tr_{\varrho}(P')=[6]\Tr_{\varrho}(P)=O$.  On the other hand, we also have $\Tr_{\varrho}(P')=\Tr_{\varrho}(P'_{1})+\Tr_{\varrho}(P'_{2})+\Tr_{\varrho}(P'_{3})=[6]P'_{1} + O + O=[6]P'_{1}$.  This implies that $[6]P'_{1}=O$.  Also, we have $\varrho^3(P')=-P'$, and thus
\begin{align*}
O=\varrho^3(P') + P' &=\varrho^3(P'_{1} + P'_{2} + P'_{3}) 
+ (P'_{1} +P'_{2} + P'_{3}) \\
&= (P'_{1} - P'_{2} + P'_{3}) + (P'_{1} + P'_{2} + P'_{3})
= [2]P'_{1} + [2]P'_{3}.
\end{align*}
This implies that $[6]P'_{3}=O$.  Similarly, since $P', \varrho^4(P')$ and $\varrho^2(P')$ are collinear, we have
\[
O= (1+\varrho^4+\varrho^2)(P')=
(1+\varrho^4+\varrho^2)(P'_{1} + P'_{2} + P'_{3})
=[3]P'_{1} + [3]P'_{2},
\]
and so $[6]P'_{2}=O$.  As a result, we have $[6]P'=[36]P=O$, and $P$ does have an infinite order.
\end{proof}
 
For a character $\chi$ of order $6$ of $\Gal(K_6/\Q)$, Lemma~\ref{lem:rankE(K6)} shows that if $P\in E(K_6)^{\chi}$ is a point of infinite order, we must  have $\rk E(K_6)^{\chi} \geq 1$.  It follows from Kato-Kolyvagin that $L(E,1,\chi)=0$ for such a character $\chi.$

\section{Explicit formulas}
Now, suppose $E$ is realized as a plane cubic curve defined by the Weierstrass equation 
\[
E:y^{2}=x^{3} + A x + B. 
\]
The origin $O$ is the point at infinity $(0:1:0)$, which is an inflection point.  The inversion map $P\mapsto -P$ on $E$ is given by $(x:y:z)\mapsto (x:-y:z)$, which is obviously the restriction of an involution $\iota:(x:y:z)\mapsto(x:-y:z)$ of~$\P^{2}$.  Let $L$ be the line given by the equation 
\[L:y=tx+uz,\] where $t,u$ are parameters.  Then the image of $L$ under the involution $\iota$ is given by $y=-tx-uz$. The degenerate conic $C_{L}=L\sqcup \iota(L)$ is thus given by \[C_{L}:(y-tx-yz)(y+tx+uz)=0.\]

If we set $\iota(u)=-u$ and $\iota(t)=-t$, it is easy to see that $[\Q(t ,u):\Q(t ,u)^{\iota}]=2$, and $\Q(t ,u)^{\iota}=\Q({t}/{u},t u)$.   For the reason explained below, we let \[T =-{u}/{t}, \quad U = -u t.\] Then the equation of $C_{L}$ is
\[
C_{L}: T y^2  - U (x - Tz)^2 = 0.
\]
The $x$-coordinate of the points of intersection $E\cap C_{L}$ satisfies the cubic equation 
\begin{equation}\label{eq:x-coord}
U (x - T)^2=T(x^3+Ax+B)
\end{equation}
By Lemma~\ref{lem:cyclic-sextic}, each point of intersection $E\cap C_{L}$ is defined over a cyclic sextic extension of~$\Q$ if and only if $x$-coordinates of points of intersection are defined over a cyclic cubic equation.  This condition is nothing but the condition that the discriminant $\Delta$ of \eqref{eq:x-coord} in $x$ is a square.  The discriminant of \eqref{eq:x-coord} is given by
\begin{multline}\label{eq:disc}
\Delta = 
4 T (T^{3}+A T +B ) U^{3}
- T^{2} (27 T^{4} + 30 A T^{2} - A^{2} + 36 B T ) U^{2}
\\
-6  T^{3} (4 A^{2} T -9 B T^{2}+3 A B ) U 
-(4 A^{3}+27 B^{2}) T^{4}.
\end{multline}
Thus, for a pair $(U,T)\in \Q\times\Q$, $\Delta$ is a square if and only if each point of $E\cap C_{L}$ is defined over some cyclic sextic extension.  
Thus, letting $\Delta/T^{2}=D^{2}$, we obtain an equation in $(U,D,T)$ in the affine space $\A^{3}/\Q$, which yields an affine model of $S_{6}(E)$ over $\Q$.
\begin{multline}\label{eq:S6-std}
S_{6}(E) :
TD^{2} 
=4 (T^{3}+A T +B ) U^{3}
- T (27 T^{4} + 30 A T^{2} - A^{2} + 36 B T ) U^{2}
\\
-6  T^{2} (4 A^{2} T -9 B T^{2}+3 A B ) U 
-(4 A^{3}+27 B^{2}) T^{3}.
\end{multline}

\begin{proposition}\label{prop:sextic-field}
Suppose $(U,D,T)$ is a $\Q$-rational point on the affine surface defined by~\eqref{eq:S6-std}. Then $E$ has a point $P=(x,y)$ such that $x$ satisfies the cubic equation 
\begin{equation}\label{eq:cubic}
T\,x^3 - U\,x^2 + T(A + 2U)\,x + T(B - TU)=0,
\end{equation}
and $y$ satisfies the sextic equation 
\begin{multline}\label{eq:sextic}
T^3 y^6 - U (3 T^4 - 2 (A + 3 U) T^2 - U^2) y^4 
\\+ U^2 (3 T^5+2 U T^3-6 B T^2+A (A+2 U) T+2 B U) y^2 
\\
- (T^3 + A T + B)^2 U^3 = 0.
\end{multline}
The splitting field $K_{3}$ of \eqref{eq:cubic} is generically a cyclic cubic extension of $\Q$, and the splitting field $K_{6}$
of \eqref{eq:sextic} is generically a cyclic sextic extension over $\Q$.
\end{proposition}

\begin{proof}
The equation \eqref{eq:S6-std} of $S_{6}(E)$ is the condition that the Galois group of \eqref{eq:cubic} is cyclic of order~$3$ as long as \eqref{eq:cubic} is irreducible over~$\Q$.  The fact that $K_{6}/\Q$ is cyclic follows from Lemma~\ref{lem:cyclic-sextic}.
\end{proof}

\begin{remark}\label{rmk:sextic}
Over the function field $\Q(T,U)$, the splitting field of the sextic equation \eqref{eq:sextic} in $y$ is a dihedral extension.  
\end{remark}

Looking at \eqref{eq:S6-std}, we see immediately that the projection $\pi_{6}:(U,D,T)\mapsto T$ defines an elliptic fibration with a section at $U=\infty$.  It turns out that there is a  $3$-torsion section defined over $\Q(\omega)$.  So, we rewrite the Weierstrass form such that the $X$-coordinate of this $3$-torsion section equals~$0$:
\begin{proposition}\label{prop:elliptic-pencil}
The projection $\pi_{6}:(U,D,T)\mapsto T$ extends to an elliptic fibration $\tilde\pi_{6}:\widetilde{S}_{6}(E)\to \P^{1}$ whose Weierstrass equation is given by
\begin{multline}\label{eq:ET-std}
\mathscr{E}_{T}:
Y^{2} = 
X^{3}
-27\bigl((3T^2 + A)X 
\\
+ 4(27 A T^4 + 54 B T^3 + 18 A^2 T^2 + 54 A B T - A^3 + 27 B^2)\bigr)^2.
\end{multline}
The singular fibers of the elliptic surface $\mathscr{E}_{T}$ are as follows.
\[
\begin{cases}
\ \mathrm{I}_{6} & \text{ at }\ T=\infty,
\\
\ \mathrm{I}_{3} & \text{ at }\ 27 A T^4 + 54 B T^3 + 18 A^2 T^2 + 54 A B T - A^3 + 27 B^2 = 0,
\\
\ \mathrm{I}_{2} & \text{ at }\ T^{3} + A T + B =0.
\end{cases}
\]
The Mordell-Weil group $\mathscr{E}_{T}\bigl(\overline{\Q(A,B)}(T)\bigr)$ contains a subgroup isomorphic to $\Z\oplus\Z/3\Z$ generated by the sections
\begin{align*}
&\mathcal{P}_{\infty} = \bigl(-12(3 A T^2 + 9 B T - A^2), 
108\sqrt{-4 A^3 - 27 B^2} (T^3 + A T + B) \bigr),
\\
&\mathcal{T}_{3} = \bigl(0, 12\sqrt{-3}\,(-27 A T^4 - 54 B T^3 - 18 A^2 T^2 - 54 A B T + A^3 - 27 B^2)\bigr).
\end{align*}
If $E$ does not have complex multiplication, the subgroup $\<\mathcal{P}_{\infty},\mathcal{T}_{3}\>$ equals $\mathscr{E}_{T}(\overline{\Q(A,B)}(T))$.
\end{proposition}

\begin{proof}
Since \eqref{eq:S6-std} is already cubic in $U$, it is straight forward to make it to the Weierstrass form.  Once this is done, we can easily compute its $3$-division polynomial, and verify that it has a root over $\Q(A,B)(T)$.  We then translate the root to $X=0$, and we have the equation~\eqref{eq:ET-std}.  Overall the change of variables is given by
\begin{equation*}
U=\frac{\bigl(X + 12(3AT^2 + 9BT - A^2)\bigr)T}{36(T^3 + AT + B)},
\quad
D=\frac{YT}{108(T^3 + AT + B)}.
\end{equation*}
The determination of the singular fibers is done by Tate's algorithm.
The section $\mathcal{T}_{3}$ is the $3$-torsion section mentioned above, 
and $\mathcal{P}_{\infty}$ of infinite order comes from the point $(U,D)=\bigl(0,\sqrt{-4 A^3 - 27 B^2}\bigr)$ in~\eqref{eq:S6-std}.  From the types and numbers of singular fibers, the Shioda-Tate formula tells us that 
\begin{align*}
\rank \mathscr{E}_{T}\bigl(\overline{\Q(A,B)}(T)\bigr)
&= \rho(\widetilde{S}_{6}) - \underbrace{2}_{\substack{\text{good fiber}\\+\text{0-section}}}
- \underbrace{(6-1)}_{\mathrm{I}_{6}} 
- \underbrace{4\times (3-1)}_{4\mathrm{I}_{3}}
- \underbrace{3\times (2-1)}_{3\mathrm{I}_{2}}
\\
&=\rho(\widetilde{S}_{6}) - 18,
\end{align*}
where $\rho(\widetilde{S}_{6})$ is the Picard number of $\widetilde{S}_{6}$.  
Thus, if $E$ does not have complex multiplication, $\rank \mathscr{E}_{T}=1$, and an easy computation of height shows that $\mathcal{P}_{\infty}$ is a generator.  The torsion subgroup of the Mordell-Weil group can be verified to be isomorphic to~$\Z/3\Z$ based on the configuration of the singular fibers.
\end{proof}

\begin{remark}\label{rmk:3-isog}
The elliptic surface $\mathcal{E}_{T}$ admits a $3$-isogeny $\mathcal{E}_{T}\to\mathcal{E}'_{T}=\mathcal{E}_{T}/\<\mathcal{T}_{3}\>$, where $\mathcal{E}'_{T}$ is given by
\[
\mathscr{E}'_{T}:
Y^2 = X^3 + \bigl((3 T^2 + A)X + 4(T^3 + AT + B)^2\bigr)^2
\]
\end{remark}

\section{Main result: General case}
\label{sec:main_result}

Based on the following observation, we can prove an unconditional result that is stronger than our cubic twist results.

\begin{lemma}\label{lem:twist}
Let $E$ be an elliptic curve over $\Q$, and $E^{\delta}$ its quadratic twist by a rational number~$\delta$.  Then the surfaces $S_{6}(E)$ and $S_{6}(E^{\delta})$ are isomorphic over~$\Q$.
\end{lemma}

\begin{proof}
If $E$ is given by $y^{2}=x^{3}+Ax+B$, then $E^{\delta}$ is defined by $\delta y^2=x^{3}+Ax+B$.  The intersection $E^{\delta}\cap C_{L}$ is given by 
\begin{equation}\label{CL-twist}
\delta U (x - T)^2=T(x^3+Ax+B).
\end{equation}
Since the equation of $S_{6}(E^{\delta})$ is obtained as the discriminant of this cubic equation, it is isomorphic to $S_{6}(E)$ via the map $\delta U\mapsto U$.
\end{proof}

In the paper \cite{FKK}, we studied the $K3$ surface defined as the minimal model of the quotient $E\times E$ divided by the automorphism 
$\rho':(P,Q)\mapsto (Q,-P-Q)$ of order~$3$.  Let $\upsilon$ be the automorphism of $E\times E$ defined by $\upsilon:(P,Q)\mapsto (P,P+Q)$.  Then we have $\rho'=\upsilon^{-1}\circ\rho^{2}\circ\upsilon$, and thus $E\times E/\<\rho'\>=E\times E/\<\rho^{2}\>$.
An affine equation of $E\times E/\<\rho'\>$ is given by the discriminant of the intersection between $E$ and the line $y=tx+u$.  
\begin{multline}\label{eq:S_3}
S_{3}(E) : 
	d^{2} = -27 u^{4} - 4 t^{3} u^{3} - 6(5 A t^{2} - 9 B) u^{2} 
	\\
	- 4 t ( A t^{4} - 9 t^{2} - 6 A^{2} ) u 
	+ 4 B t^{6} + A^{2} t^{4} - 18 A B t^{2} - ( 4 A^{3} + 27 B^{2} ).
\end{multline}
We have proved (\cite[\S5]{FKK}) that if $E(\Q)$ contains at least six points, then  $S_{3}(E)(\Q)$ is Zariski dense.

By construction, $S_{3}(E)$ is a double cover of $S_{6}(E)$, and the covering map $\pi_{2}:S_{3}(E) \to S_{6}(E)$ is given by $\pi_{2}:(t,u,d)\mapsto (T,U,D)=(-u/t, -u t, -d u/t)$.  Unlike $S_{6}(E)$, $S_{3}(E^{\delta})$ depends on the value of~$\delta$.

\begin{theorem}\label{thm:main}
Let $E$ be an elliptic curve defined over $\Q$, and let $\delta$ be a square-free integer such that $\rank E^{\delta}(\Q)\ge 1$.
Then there are infinitely many cubic extensions $K_{6}/\Q(\sqrt{\delta})$ such that $K_{6}/\Q$ is cyclic sextic, and $\rank E(K_{6})$ is strictly greater than the rank of the subgroup of $E(K_6)$ generated by $E(F)$ for all proper subfields $F\subset K_6$.
\end{theorem}

\begin{proof}
By \cite[Th.~5.1]{FKK}, $S_{3}(E^{\delta})(\Q)$ contains a rational point $(t,u,d)$ that corresponds to a 
 point $P=(x_{P},y_{P})\in E^{\delta}(K_{3})\setminus E^{\delta}(\Q)$ where $K_{3}$ is a cyclic cubic extension.  Then, the point $(t,u,d)\in S_{3}(E^{\delta})(\Q)$ maps to $(T,U,D)=(-u/t, -u t,-d u/t)\in S_{6}(E^{\delta})(\Q)$.  It corresponds to a point $(T,U/\delta,D)=(-{u}/{t}, -u t/\delta, -d u/t)$ in $S_{6}(E)(\Q)$ via the correspondence $\delta U\mapsto U$ in Lemma~\ref{lem:twist}.  This implies that $K_{6}=K_{3}(\sqrt{\delta})$ is a cyclic sextic extension, and the corresponding point $(x_{P},y_{P}/\sqrt{\delta})$ on $E$ is defined over~$K_{6}$. 

Since there are infinitely many different cyclic cubic extensions $K_{3}$ such that $\rank E^{\delta}(K_{3})$ is greater than $\rank E^{\delta}(\Q)$ (\cite[Th.~5.1]{FKK}),  there are as many different cyclic sextic extensions $K_{6}$ with the desired property for each such $\delta$.  
\end{proof}

\begin{corollary}\label{cor:main}
For any elliptic curve over $\Q$, there are infinitely many cyclic sextic extensions $K_{6}$ such that $\rank E(K_{6})$ is strictly greater than the rank of the subgroup of $E(K_6)$ generated by $E(F)$ for all proper subfields $F\subset K_6.$
\end{corollary}

\begin{proof}
It is easy to see that there are infinitely many distinct square-free integers $\delta$ such that $\rank E^{\delta}(\Q)$ is positive.
\end{proof}

\begin{corollary}\label{cor:Dirichlet-char}
For any elliptic curve $E/\Q$, there are infinitely many primitive Dirichlet characters $\chi$ of order $6$ such that $L(E,1,\chi)=0.$
\end{corollary}

This answers a question posed in \cite{Berg-Ryan-Young}.

\section{Various special cases}

Since the results in \S\ref{sec:main_result} inherits the drawback of the results in \cite[\S5]{FKK}, whose proofs heavily rely on Faltings's theorem and 
Merel’s theorem on the bound for the torsion points of an elliptic curve, they do not yield quantitative results.  In this section, we look at some special cases and count the number of different sextic fields over which $E$ has a positive rank.
 
\subsection{$E$ with a rational $3$-isogeny}

Suppose $E$ has a rational $3$-isogeny.  If so, $E$ can be written in the form 
\begin{equation}\label{eq:E-3-tors}
E_{a,b}: y^{2} = x^3 + a(x-b)^2, \quad a, b\neq 0,
\end{equation}
with $3$-torsion points $(0,\pm b\sqrt{a})$.  In this case the equation of $S_{6}(E_{a,b})$ is given by
\begin{multline}\label{eq:S6-3-isog}
S_{6}(E_{a,b}) :
T D^{2} = 
4 \left(T^{3}+ a(T - b)^{2}\right) U^{3}
\\
-T\left(27 T^{4}+36 a T^{3} +4 a (2 a -15 b ) T^{2}-4 a b (4 a -9 b ) T +8 a^{2} b^{2}\right) U^{2}
\\
+2 a T^{2}\left((2 a^{2}+18 a b +27 b^{2}) T^{2}-2 a b (2 a +15 b ) T 
+2 a b^{2} (a +9 b )\right) U 
\\
-a^{2} b^{3} (4 a +27 b )\, T^{3}.
\end{multline}
The discriminant of the right hand side as a cubic polynomial in $U$ is 
\[
16 a^{3} T^{6} (T - b)^{3}\bigl(54 b T^3+ 8a^{2}(T - b)^2 + 9 a (T - b)(T^2 + 6 b T - 3 b^2)\bigr)^{3}.
\]
This means that the fiber at $T=b$ of the elliptic fibration $\mathscr{E}_{T}$ is singular.  It is easy to see that its Kodaira type is $\mathrm{I}_{3}$.  
If we let $T=b$ in \eqref{eq:S6-3-isog}, we have
\[
b D^2 = b^{3}(4 U -4 a b - 27 b^2)(U -a b )^2.
\]
This means that the fiber at $T=b$ in the affine model \eqref{eq:S6-3-isog} is a parametrizable curve, and we can easily find its parametrization:
\[
s\mapsto (U,D,T) = \left(s^{2} + a b + \tfrac{27}{4} b^2, 2b s^3 + \tfrac{27}{2} b^3 s,b\right)\in S_{6}(E_{a,b})(\Qbar).
\]
Let $[P,Q]\in E_{a,b}\times E_{a,b}/\<\rho\>$ be the point corresponding to the image of~$s$.  Then, the $x$ coordinates of $P$ and $Q$ are roots of the cubic equation
\begin{equation}\label{eq:cubic-3-tors}
4 x^{3}- b (s^{2}+27) x^{2} +2 b^{2}  (s^{2}+27) x - b^{3}(s^{2}+27) =0,
\end{equation}
and the $y$ coordinates satisfy the quadratic equation
\[
4 y^{2} =(b s^{2}+4 a +27 b)(x -b)^{2}.
\]
Let $s=3m/n$, and $x=3bx'$ in \eqref{eq:cubic-3-tors}, and define 
\[
f_{a,b}(x')= 12 n^2 x'^3 - 9 (m^2 + 3 n^2)x'^2 
+ 6 (m^2 + 3n^2) x' - (m^2 + 3n^2).
\]
It is irreducible over $\Q(m,n)$ as $f_{a,b}(x')$ is Eisenstein at $(m^2 + 3n^{2})$, and has discriminant
\[
\operatorname{Disc}(f_{a,b}) = 2^4 3^4 n^2 m^2 (m^2 + 3n^2)^2.
\]
Since $\operatorname{Disc}(f_{a,b})$ is a square, the splitting field $K_{3,m,n}/\Q(m,n)$ of $f_{a,b}(x')$ is generically a cyclic cubic extension. 

We now specialize $m$ and $n$  to be coprime integers such that $(m^2 + 3n^2)>1$ is square-free and consider the specialized extension $K_{3,m,n}/\Q$ (using the same notation). Then $f_{a,b}(x')$ is Eisenstein at any prime $p$ dividing 
 $(m^2 + 3n^2)$, so that  $K_{3,m,n}/\Q$  is a cyclic cubic extension of $\Q.$ 
For any prime $p$ dividing $mn$, $f_{a,b}(x')\bmod p$ has a double root at $x' = 1$, so that there is a prime of $K_{3,m,n}$ over $p$ which has degree~$1$, and is therefore unramified in $K_{3,m,n}$.  Consequently the Discriminant  $D(K_{3,m,n}) = (m^2 +3n^{2})^2$ (up to constants depending only on $E_{a,b}$), and so the conductor of $K_{3,m,n}$ is $m^2 + 3n^{2}$ (up to constants).

Also the $y$-coordinate lies in $K_{6,a,b,m,n}=K_{3,a,b}\bigl(\sqrt{9 b m^2 + (4 a + 27 b)n^{2}}\bigr)$.  If $9 b m^2 + (4 a + 27 b)n^{2}$ is square-free, then $K_{6,a,b,m,n}$ is a quadratic extension of $K_{3,m,n}$ of conductor $9 b m^2 + (4 a + 27 b)n^{2}$. 

We note that the only common prime factors of $(9bm^2 + (4a+27b)n^2)$ and $(m^2 +3n^2)$ must be divisors of $6ab$. This follows since  $(9bm^2 + (4a+27b)n^2) - 9b(m^2 +3n^2) = 4an^2$ and so if $p$ is a prime dividing both factors, then $p$ divides $4an$. If in addition $p$ divides $n$ then  $p$ also divides $9bm^4$ (reading modulo $p$).
So finally, the sextic field $K_{6,a,b,m,n}$ is a cyclic extension of $\Q$ with conductor 
\[
g_{a,b}(m,n) =  (9 b m^2 + (4 a + 27 b)n^{2})(m^2 + 3n^{2})
\]
up to factors of $6ab$. 
Thus, the number of such fields (or pairs of sextic characters) of conductor less than~$X$ is at least as large as the number of distinct square-free conductors 
less than~$X$ of cyclic sextic fields, which is in turn at least equal to (up to a possible factor dividing $6ab$) the number of square-free values less than $X$ taken by the form $g_{a,b}(m,n)$.  As in the case of cyclic cubics, we use the following theorem by Stewart and Top \cite{St-T} to give an estimate.  For convenience, we quote it here in its full form with the original notation.

\begin{theorem}[{Stewart-Top\cite[Theorem~1]{St-T}}]\label{thm:Stewart-Top}
Let $A, B, M$ and $k$ be integers with $M \geq 1$ and
$k \geq 2$.  Let
\begin{equation}
F(X, Y) = a_{r} X^{r} + a_{r} X^{r - 1} Y + \cdots + a_{0} Y^{r}
\end{equation}
be a binary form with integer coefficients,
non-zero discriminant and degree $r$ with $r \geq 3.$  Let $m$ be the
largest degree of an irreducible factor of $F$ over $\Q$ and
suppose that $m \leq 2k + 1$ or that $k = 2$ and $m = 6$. Let $w$ be the largest positive integer such that $w^k$ divides $F(h_1,h_2)$ for all integers $h_1$ and $h_2$ with  $h_1 \equiv A \Mod M$ and $\ h_2 \equiv B \Mod M.$
 For any real number $x$ let $R_{k}(x)$ denote the number of $k$-free
integers $t$ with $|t| \leq x$ for which there are integers $h_1$ and $h_2$, 
  with 
\[ h_1 \equiv A \Mod M, \ h_2 \equiv B \Mod M,
 \text{ and } F(h_1, h_2) = t w^{k}.
 \]
Then, there are positive
numbers $C_{1}$ and $C_{2}$ which depend on $M, k$ and $F,$ such that
if $x$ is a real number larger than $C_{1}$, then
\[
R_{k} (x) > C_{2}~ x^{\frac{2}{r}}.
\]
\end{theorem}

\begin{theorem}\label{thm:3-isog}
Let $E_{a,b}$ be an elliptic curve with an equation of the form \eqref{eq:E-3-tors}. 
Then the number of sextic Dirichlet characters $\chi$ of conductor less than $X$ for which $L(E_{a,b},1,\chi)=0$ is $\gg X^{\frac{1}{2}}$.
\end{theorem}

\begin{proof}
In Theorem~\ref{thm:Stewart-Top},  use $g_{a,b}(m,n)$ as $F(X,Y)$ and let $M=1$ and $k=2$.  Taking $(h_1,h_2) = (0,1)$ and $(1,0)$ we see that the integer $w$ in Theorem~\ref{thm:Stewart-Top} must divide $\gcd (9b,12a)$ and so is bounded for  $E_{a,b}.$
Then, its assertion implies the number of distinct square-free values less than $X$ taken by the form $g_{a,b}(m,n)$ is  $\gg X^{\frac{1}{2}}$.  Thus, the number of sextic fields having desired conditions is at least as large.
\end{proof}

\begin{corollary}\label{cor:3-isog}
The number $V_{\Psi_6,E_{a,b}}(X)$ of sextic fields $K_{6}$ of conductor less than $X$ for which $\rank E_{a,b}(K_{6})$ is strictly greater than the rank of the subgroup of $E_{a,b}(K_6)$ generated by $E(F)$  for all proper subfields $F\subset K_6$  satisfies $V_{\Psi_6,E_{a,b}}(X) \gg X^{\frac{1}{2}}$.
\end{corollary}

\subsection{$E$ with a rational $2$-torsion point}

Suppose $E$ has a rational $2$-torsion point.  If so, $E$ can be written in the form 
\[
E: y^{2} = x(x^2 + a x + b),
\]
with $2$-torsion point $(0,0)$.  In this case the equation of $S_{6}$ is given by
\begin{multline}\label{eq:S6-2-torsion-T}
S_{6}(E): D^{2} = 
4 \left(T^{2}+ a T +b \right) U^{3}
\\
-\left(27 T^{4}+36 a T^{3} + 2(4 a^{2}+15 b) T^{2} +8  a b T -b^{2}\right) U^{2}
\\
+2 T \left(a(2a^2 - 9b)T^{2} + 2b(a^2 - 6b)T - a b^2 \right) T^{2},
\end{multline}
which can be viewed as an elliptic curve over $\Q(T)$.  Collecting terms with respect to $T$, we have
\begin{multline}\label{eq:S6-2-torsion-U}
D^{2} = 
-27 U^{2} T^{4} -2 a U \left(18 U - 2 a^{2} +9 b \right) T^{3}
\\
+\left(4 U^{3} - 2(4a^2 + 15b)U^{2} + 4b(a^2 - 6b)U + b^2(a^2 - 4b)\right) T^{2} 
\\
+2 U \left(2 U^{2}-4 b U -b^{2}\right) a T 
+ b U^{2} \left(4 U +b \right),
\end{multline}
which can be viewed as a curve of genus~$1$ over $\Q(U)$.  Calculating the discriminant of the right hand side of \eqref{eq:S6-2-torsion-U} with respect to~$T$, we have
\begin{multline*}
\Delta=
256U^4(a^2 - 4b)(16 U -a^2 + 4b)
\\
\times\left(3U^{4} - 8(a^2 - 3b)U^{3} 
- 6b(2a^2 - 9b)U^2 - 6b^2(a^2 - 4b)U - b^3(a^2 - 3b)\right)^3.
\end{multline*}
This shows that the fiber of $(T,U,D)\mapsto U$ at $U=(a^{2}-4b)/16$ is singular.  Indeed, it is given by
\begin{equation}\label{eq:U-fiber}
2^{10} D^{2} = 
- (a^{2}-4 b)(2 T -a )^{2} \left(27(a^2 - 4b)T^2 - a(a^2 - 36 b)T - b(a^2 - 4 b)\right).
\end{equation}
It may be rewritten in the following form:
\[
\left(54(a^2 - 4b)T-a(a^2 - 36b)\right)^{2}
+ 3\left(2^{6}3 D\right)^{2}
=(a^{2}+12b)^{3}.
\]
Thus, if $a^{2}+12b$ is a norm in $\Q(\sqrt{-3})$, then \eqref{eq:U-fiber} can be parametrized over $\Q$.  Instead of working with the general case, we work out in detail with an example.

Let $a=-4$ and $b=-1$.  The curve 
\[
E_{\mathrm{160b1}}:y^{2}=x(x^{2}-4x-1)
\]
is a curve of conductor $160$, and labeled ``160b1'' in the Cremona list.  
Its Mordell-Weil group equals $\Z/2\Z$.  
The surface $S_{6}(E_{\mathrm{160b1}})$ is given by
\begin{multline*}
S_{6}(E_{\mathrm{160b1}}): D^{2}=-27 T^4 U^2 + 8U(18U - 41)T^3 
\\
+ 2(2U^3 - 49U^2 - 44U + 10)T^2 - 8U(2U^2 + 4U - 1)T - U^2(4U - 1),
\end{multline*}
and the fiber at $U=5/4$ is singular and given by 
\[
16D^{2} = -5(5T + 1)(27T + 5)(T + 2)^2.
\]
It is now straightforward to parametrize this curve, and we see
that $S_{6}(E_{\mathrm{160b1}})$ contains a parametrized curve 
\begin{equation}\label{eq:E160b1-param}
m/n\mapsto (U,D,T)=\left(\frac{5}{4},
\frac{mn(245m^2 + 81 n^{2})}{270(3m^2 + n^{2})^2},
-\frac{25m^2 + 9 n^{2}}{45(3m^2 + n^{2})}\right) \in S_{6}(E_{\mathrm{160b1}}).
\end{equation}
From this we see that for each value of $m$, there is a point $P(x,y)$ on $E_{\mathrm{160b1}}$ whose $x$ coordinate satisfies the equation
\begin{multline*}
f_{m/n}(x)=36 (3 m^2 + 25 n^{2}) (m^2 + 9 n^{2}) x^3  
+ 9 (3 m^2 + 25 n^{2}) (11 m^2 + 81 n^{2}) x^2  
\\
+ 54 (3 m^2 + 25 n^{2}) (m^2 + 9 n^{2}) x  + (m^2 + 9 n^{2})^2=0.
\end{multline*}
Let $K_{3}$ be the splitting field of $f_{m/n}(x)$. Since the discriminant of $f_{m/n}(x)$ is a square, $K_{3}/\Q$ is a cyclic cubic extension, and $K_{3}=\Q[x]/(f(x))$.  

\begin{lemma}\label{lem:K3-ramif}
Let $p$ be a prime number different from $2,3$ or $5$. 
\begin{itemize}
\item 
If $p$ divides $3m^2 + 25n^2$ square-freely, then $p$ ramifies (totally) in $K_3=\Q[x]/(f_{m/n}(x))$.  
\item 
If $p$ is prime to $3m^2 + 25n^2$ then $p$ is unramified in $K_3$.  
\item
If $3m^2 + 25n^2$ is square-free, then $K_3/\Q$ is a cyclic cubic extension of conductor $3m^2 + 25n^2$ up to factors of $180$. 
\end{itemize}
\end{lemma}

\begin{proof}
We first note that the polynomial $f_{m/n}(x)$ has discriminant
\[
\operatorname{Disc}(f_{m/n})
=2^{6}3^{6}m^2 n^{2}(49m^2 + 405 n^{2})^2
(3m^2 + 25 n^{2})^2 (m^2 + 9 n^{2})^2.
\]
and is ``reverse Eisenstein''
for $3m^2 + 25n^2$. Hence any prime $p$ dividing $3m^2 + 25n^2$ square-freely is ramified in $K_3$. 
 
 Now consider the change of variables $x_{1}=36 (m^2 + 9 n^2)(3 m^2 + 25 n^2)x$. Then $f_{m/n}(x)$ becomes 
\begin{multline*}
f_{m/n,1}(x_{1})= x_{1}^3 + 3 (3 m^2 + 25 n^2) (11 m^2 + 81 n^2) x_{1}^2 
\\
+ 216 (3 m^2 + 25 n^2)^2 (m^2 + 9 n^2)^2 x_{1} 
+ 1200 (3 m^2 + 25 n^2)^2 (m^2 + 9 n^2)^4. 
\end{multline*}
Suppose $m$ and $n$ are relatively prime and define $F = 3 m^2 + 25 n^2$, and $G = m^2 + 9 n^2$.  Then, we have $11 m^2 + 81 n^2 = 9F - 16G$, and $49m^2 + 405 n^{2}=18F-5G$.  Moreover, we have
\[
f_{m/n,1}(x_{1})=x_{1}^{3} + 27 F^2 x_{1}^2 - 1200 F G x_{1}^2 
+ 216 F^2 G^2 x_{1} + 1200 F^2 G^4.
\]
Suppose $p$ is a prime different from $2$, $3$, and $5$ that divides the discriminant $\operatorname{Disc}(f_{m/n})=2^{6}3^{6}m^{2}n^{2}F^{2}G^{2}(18F-5G)^{2}$. If $p$ divides $F$ square-freely, we have seen that  $p$ ramifies in $K_{3}/\Q$. If $p$ divides $G$, we have $f_{m/n,1}\equiv x_{1}^{2}(x_{1}+27F^{2})\pmod p$. Since $F-3G=2n^{2}$, if $p$ divides also $F$, it divides both $m$ and $n$, which is a contradiction.  Thus, $p$ does not divide $F$, and $f_{m/n,1}$ has a degree 1 factor modulo $p$. Since the extension $K_3$ is a cyclic cubic extension, $p$ splits in $K_3$ and must be unramified.
If $p$ divides $m$, then we have $f_{m/n,1}(x_{1})\equiv (x_{1} + 2700 n^4)^2(x_{1} + 675 n^4)\pmod p$.  
If $p\neq 3, 5$, then $f_{m/n,1}$ has a $\text{degree~$1$ factor}\bmod p$, and $p$ is unramified in $K_3$.  If $p$ divides $n$, then we have $f_{m/n,1}(x_{1})\equiv (x_{1} + 12 m^4)^2(x_{1} + 75 m^4)$.  So, if $p\neq 3$, $p$ is unramified.  Finally, if $p$ divide $49m^2 + 405 n^{2}=18F-5G$, then $f_{m/n,1}(x_{1})\equiv 1/25(5 x_{1} - 432 F^2)^2(x_{1} + 27 F^2) \pmod p$.  If $p\neq 2, 3, 5$, then $f_{m/n,1}$ has a degree 1 factor modulo $p$, and $p$ is unramified.  
\end{proof}

The $y$ coordinate of $P(x,y)$ satisfies
\[
36 y^2 (m^2 + 9 n^{2}) (3 m^2 + 25 n^{2}) - \left(9 (3 m^2 +25 n^{2}) x + 5 (m^2 + 9 n^{2})\right)^2=0.
\]
In other words, $P=(x,y)$ is defined over $K_{6}=K_{3}\bigl(\sqrt{ (m^2 + 9 n^{2}) (3 m^2 + 25 n^{2})}\bigr)$ which is a cyclic sextic extension of $\Q$.  Once again, by Stewart-Top \cite{St-T} we have the following.

\begin{proposition}\label{prop:E160b1}
The number $V_{\Psi_6,E_{\mathrm{160b1}}}(X)$ of cyclic sextic fields $K_{6}$ of conductor less than $X$ for which $\rank E_{\mathrm{160b1}}(K_{6})$ is  strictly greater than the rank of the subgroup of $E(K_6)$ generated by $E(F)$  for all proper subfields $F\subset K_6$ satisfies $V_{\Psi_6,E_{\mathrm{160b1}}}(X)\gg X^{\frac{1}{2}}$.\end{proposition}

\begin{example}\label{eg:160b1}
(1) \ 
Let $m=3, n=1$ in \eqref{eq:E160b1-param}.  We have 
\[
m^2 + 9n^2=2 \cdot 3^2, \quad 3m^2 + 25n^2=2^{2}\cdot 13,
\]
and there is a point $P(x,y)\in E_{160b1}(\Qbar)$ whose $x$-coordinate satisfies the cubic equation
\[
(2^3 \cdot 13) x^{3} + (2^2\cdot 5\cdot 13) x^{2} + (2^2\cdot 3\cdot 13) x + 5^2 = 0,
\]
and its $y$-coordinate is given by 
\(
y =\frac{\sqrt{-26}}{52}(26 x + 5).
\)
The discriminant of the cubic equation above equals $2^{6}\cdot 13^{2}\cdot 47^{2}$.
In this case, it turns out that $K_{3}=\Q(\zeta_{13}+\zeta_{13}^{-1}+\zeta_{13}^5+\zeta_{13}^{-5})=\Q[x]/(x^3 + x^2 - 4x + 1)$, where $\zeta_{13}$ is a primitive 13th root of unity.  
In fact, one of the points $P(x,y)$ is given by
\[
P(x,y)=\biggl(-\frac{(\alpha - 4)^2}{26}, 
\frac{\bigl((\alpha - 4)^2 - 5\bigr)\sqrt{-26}}{52}\biggr),
\text{ where $\alpha=\zeta_{13}+\zeta_{13}^{-1}+\zeta_{13}^5+\zeta_{13}^{-5}$.}
\]
So, $K_{3}$ ramifies only at $13$, and $K_{6}=K_{3}(\sqrt{-26})$.  $K_{6}$ is a cyclic sextic extension ramifying at $2$, and~$13$.
\smallskip

(2) \ 
Let $m=5, n=1$ in \eqref{eq:E160b1-param}.  We have 
\[
m^2 + 9n^2=2 \cdot 17, \quad 3m^2 + 25n^2=2^{2}\cdot 5^{2},
\]
and there is a point $P(x,y)\in E_{160b1}(\Qbar)$ whose $x$-coordinate satisfies the cubic equation
\[
(2^3 \cdot 3^{2}\cdot 17) x^{3} + (2^2\cdot 3^{2}\cdot 89) x^{2} 
+ (2^2\cdot 3^{3}\cdot 17) x + 17^2 = 0,
\]
and its $y$-coordinate is given by 
\(
y =\frac{\sqrt{-34}}{204}(90 x + 17).
\)
The discriminant of the cubic equation above equals $2^{6}\cdot 3^{6} \cdot 17^{2}\cdot 163^{2}$.
In this case, it turns out that $K_{3}=\Q(\zeta_{9}+\zeta_{9}^{-1})=\Q[x]/(x^3 - 3x + 1)$, where $\zeta_{9}$ is a primitive 9th root of unity.  
$P(x,y)$ is given by
\[
P(x,y)=\biggl(-\frac{17(\beta - 8)}{6(15\beta + 43)}, 
\Bigl(\frac{5(\beta - 8)}{4(15\beta + 43)} + \frac{1}{12}\Bigr)\sqrt{-34}\biggr),
\text{ where $\beta=\zeta_{9}+\zeta_{9}^{-1}$.}
\]
Thus, $K_{3}$ ramifies only at $3$, and $K_{6}=K_{3}(\sqrt{-34})$ is a cyclic sextic extension ramifying at $2$, $3$, and~$17$.
\end{example}

\subsection{$E$ with $\Q(E[2])/\Q$ cyclic cubic}

If $E$ is given by $y^2=x^3+A x + B$, the field of definition of $E[2]$ is a cyclic cubic extension if and only if the discriminant of the cubic $x^3+A x + B$ is a square; i.e.,
\begin{equation}\label{eq:discrim}
-(4A^3 + 27B^2) = \delta^2 \quad \text{ for some $\delta\in \Q$.}
\end{equation}
Write $B = 2 b A$ and $\delta = 18 b c A$ with new parameters $b$ and~$c$.  Then, the condition \eqref{eq:discrim} becomes $A = -27 b^2(3c^2+1)$.  If this condition is satisfied, the equation of $E$ may be written in the form
\[
E_{c}^{b}: 3 b y^2 = x^3 - 3 (3c^2 + 1) b^2 x - 54 (3c^2 + 1) b^3.
\]
Starting with this model, we obtain the following equation of $S_{6}(E_{c}^{b})$.
\begin{multline*}
S_{6}(E_{c}^{b}):
T D^{2} = 
12 b^{3} (T^{3}-9 T \,c^{2}-6 c^{2}-3 T -2) U^{3}
\\
-9 T \,b^{2} (3 T^{4}-30 T^{2} c^{2}-9 c^{4}-24 T \,c^{2}-10 T^{2}-6 c^{2}-8 T -1) U^{2}
\\
-36 T^{2} (3 c^{2}+1) (6 T \,c^{2}+T^{2}+3 c^{2}+2 T +1) b U 
+36 c^{2} (3 c^{2}+1)^{2} T^{3}.
\end{multline*}
Then, the Weierstrass form of the elliptic fibration $\pi_{6}:(U,D,T)\mapsto T$ is given by
\begin{multline*}
\mathscr{E}_{c,T}:
Y^{2} = 
X^{3}
-27 \Bigl( (T^{2}- (3 c^{2}+1))X 
\\
-4 (3 c^{2}+1) \bigl(3 T^{4}+4 T^{3}-6 (3 c^{2}+1) T^{2}
\\
-12 (3 c^{2}+1) T -(3 c^{2}+5) (3 c^{2}+1)\bigr)\Bigr)^{2}
\end{multline*}
Note that $\mathscr{E}_{c,T}$ does not depend on~$b$.  This reflects the fact that $S_{6}(E_{c}^{b})\simeq S_{6}(E_{c})$ via the map $U\mapsto bU$.

Looking at the equation of $S_{6}(E_{c}^{b})$, we see that the fibration $\pi_{6}$ has an obvious section $T\mapsto (U,D,T)=\bigl(0,6c(3c^2+1)T,T\bigr)$.  This corresponds to the section of $\mathscr{E}_{c,T}$ given by
\[
P=\Bigl(12 (3 c^{2}+1)(T^{2}+ 2 T +3 c^{2} +1), 
72 c (3 c^{2}+1) (T^{3}-3 (3 c^{2}+1)T -2(3 c^{2}+1)) 
\Bigr).
\]
The curve on $S_{6}(E_{c}^{b})$ corresponding to the section $2P$ of $\mathscr{E}_{c,T}$ is
\begin{equation}\label{eq:param-curve-c}
T\mapsto (U,D,T)=\left(
\frac{3(3c^2 + 1)T(T^{2} + 2T - c^{2} + 1)(T^2 + 2T + 3c^2 + 1)}
{4bc^{2}(T^3 - 3(3c^2 + 1)T - 2(3c^2 + 1))},d_{b,c}(T),T\right),
\end{equation}
where $d_{b,c}(T)$ is a rational function in $T$.
For almost all values of $T$, the point $(U,D,T)\in S_{6}(E_{c}^{b})(\Q)$ gives a point on $E^{b}$ defined over some sextic cyclic sextic extension of $\Q$.  
The $x$-coordinate the point of $E_{c}^{b}$ corresponding to $(U,D,T)$ is a root of the cubic equation 
\[
Tx^3 - 3 b U x^2 + 3 T (2 b U - 3 c^2 - 1) x - T(3 b U T + 6c^2 + 2) = 0.
\]
Restricting this to the parametric curve \eqref{eq:param-curve-c}, and changing the parameter by $T=n/m+1$, we have
\begin{multline*}
4 (3 m^3 - 9c m^2 n - 3 m n^2 + c n^3)x^{3}
+ 9 (3 c^2 + 1) (3 m^2 + n^2) (m^{2} - n^{2}) x^{2}
\\
+ 6 (3 c^{2}+1)(3 m^5 + 9 c m^4 n + 4 c m^2 n^3 - 3 m n^4 + 3 c n^5) x 
\\
+ (3 c^{2} + 1) \Bigl(
3 m^6 + 18 c m^5 n + 3 (9 c^2 + 2)m^4 n^2 + 28 c m^3 n^3
\\
 - 9 (2 c^2 + 1) m^2 n^4 + 18 c m n^5 - 9 c^2 n^6\Bigr)=0.
\end{multline*}
 Its discriminant is
\begin{multline*}
2^{4}3^{4} (3 c^{2}+1)^{2}(3 m^3 - 9c m^2 n - 3 m n^2 + c n^3)^{2} 
\\
\times(3 m^6 + 18 c m^5 n - 3 m^4 n^2 + 28 c m^3 n^3 - 3 m^2 n^4 
+ 18 c m n^5 + 3 n^6)^{2},
\end{multline*}
and so the $x$-coordinate is defined over a cyclic cubic extension.  The $y$-coordinate satisfies
\[
y^{2}=\frac{- 3 (3 c^{2}+1)(3m^{2}+n^{2})(m^{2} - n^{2})}
{4 b m^{3} (3 m^3 - 9c m^2 n - 3 m n^2 + c n^3)}
(m x + m - c n)^{2},
\]
and it is defined over a cyclic sextic extension.

\def\arXiv#1{arXiv:\href{http://arXiv.org/abs/#1}{#1}}

\providecommand{\bysame}{\leavevmode\hbox to3em{\hrulefill}\thinspace}
\providecommand{\MR}{\relax\ifhmode\unskip\space\fi MR }
\providecommand{\MRhref}[2]{%
  \href{http://www.ams.org/mathscinet-getitem?mr=#1}{#2}
}
\providecommand{\href}[2]{#2}

\end{document}